\documentclass{article}
\usepackage[T1]{fontenc}
\usepackage{lmodern}
\usepackage{chngcntr}
\usepackage{ae,aecompl}
\usepackage{amsmath,amsthm,amssymb,enumerate,mathrsfs}
\usepackage{here,float}
\usepackage{amsfonts}
\usepackage{graphicx}
\usepackage{bbm}
\usepackage{caption}
\usepackage{url}
\usepackage{enumerate}
 \usepackage{bbm}
\usepackage{enumitem}
\usepackage{appendix}
\usepackage[latin1]{inputenc}
\usepackage[round]{natbib}
\usepackage[english]{babel}
\usepackage[a4paper]{geometry}
\usepackage{color}
\usepackage{setspace}
\usepackage{verbatim}
\usepackage{authblk}

\usepackage{xspace}
\usepackage{color}
\normalsize

\newtheorem{thm}{Theorem}

\theoremstyle{definition}
\usepackage{titling}

\theoremstyle{remark}

\theoremstyle{remark}
\usepackage[latin1]{inputenc}

\def\bkR{{\rm I\kern-.17em R}}

\def \1n{1\hskip -3pt \mbox{N}}

\begin{document}
\selectlanguage{english}

\title{Estimation of the intercept parameter in integrated Galton-Watson processes}

\author{Yang Lu\thanks{Department of Mathematics and Statistics,  Concordia University,  Montreal, Canada.  Email: yang.lu@concordia.ca}}

 \maketitle

\begin{abstract} We study estimation 
	of the intercept parameter in an integrated Galton-Watson process, a basic building-block for many count-valued time series models. In this unit root setting, the ordinary least squares estimator is inconsistent, whereas an existing weighted least squares (WLS) estimator is consistent only in the case where the process is transient, a condition that depends on the unknown intercept parameter . We propose an alternative WLS estimator based on 
	the new weight function of $1/t$, and show that it is consistent regardless of whether the process is transient or null recurrent, with a convergence rate of $\sqrt{\ln n}$.
\end{abstract}




\textbf{Acknowledgment}: We thank NSERC for grant RGPIN-2021-04144. 

\textbf{Data Availability Statement}: No data is generated. 

\textbf{Conflict of Interest Statement: } The author declares there is no conflict of interest. 

\section{Introduction}
Recently, there is a surge of interest in non-stationary integer-valued processes \citep{barczy2011asymptotic, barczy2014asymptotic, michel2020limiting, pei2023forecasting}. 
  Many of these models can be written as higher-order extensions of the  Galton-Watson process or branching process with immigration:
\begin{equation}
\label{gw}
X_0=0, \qquad   X_t= \sum_{i=1}^{X_{t-1}}Z_{i,t}+\epsilon_t, \forall t>0,
\end{equation}
where $Z_{i,t}, i, t$ varying, are i.i.d. with mean $m$ and variance $\sigma^2>0$, and $\epsilon_t$, $t$ varying, are i.i.d. with mean $\mu>0$ and variance $b>0$,  and are independent from  the sequence $(Z_{i,t})$. This model leads to:
$$
\mathbb{E}[X_t \mid X_{t-1}]=m X_{t-1}+\mu, 
$$
and we say that process $(X_t)$ has a unit root, or is integrated, if $m=1$. This family of integrated processes includes the INARCH model (with $Z_{i,t}$ and $\epsilon_t$ are both Poisson distributed), the NBAR model of \cite{gourieroux2019negative} (with $Z_{i,t}$ geometrically distributed and $\epsilon_t$ negative binomial distributed, sharing the same probability parameter) and many others [see e.g.  \cite{lu2021predictive} and the references therein].  However, under the unit root condition $m=1$, the INAR(1) model is excluded. Indeed, when $m=1$, we have $Z_{i,t}=1$ almost surely. Thus its variance $\sigma^2$ is equal to zero, which is disallowed. In this case, the INAR(1) model can be written as: $X_t=X_{t-1}+\epsilon_t$.  Then $\mu=\mathbb{E}[\epsilon_t]$ can be trivially estimated by $\frac{X_n}{n}$ at the standard parametric rate $\sqrt{n}$. 

One important gap in this literature is that up to now, for many of these nonstationary models, the Ordinary Least Squares (OLS) estimator is not consistent. In particular, it is known since at least \cite{winnicki1986unified} that the OLS estimator of $(m,\mu)$ of the baseline Markov model \eqref{gw}, defined by:
\begin{equation}
\label{cls}
(\hat{m}, \hat{\mu})':=\arg\min_{(m, \mu)'} \sum_{t=1}^n (X_t-m X_{t-1}-\mu)^2,
\end{equation}
is such that $\hat{\mu}$ is inconsistent (see Theorem 3 below for its large sample property).  In an attempt to solve this inconsistency issue, \cite{wei1990estimation} consider the following WLS estimator:
\begin{align}
\label{1stwls}
(\hat{\hat{m}},\hat{\hat{\mu}})':=\arg \min_{(m, \mu)} \sum_{t=1}^n \frac{(X_t-m X_{t-1}-\mu)^2}{1+X_{t-1}}. 
\end{align}
They are only able to show [see their Theorem 2.5  and \cite[Theorem 2.1]{wei1991convergence}] the consistency of $\hat{\hat{\mu}}$ (at the convergence rate of $\sqrt{\ln n}$), under the additional condition on the parameters:
\begin{equation}
\label{transient}
\tau_0:=\frac{2\mu_0}{\sigma^2_0}>1,
\end{equation}
where $\mu_0$ (resp.  $\sigma_0$) denotes the true parameter value of $\mu$ (resp.  $\sigma$).  This inequality can be interpreted \citep{wei1991convergence} as the transience condition of the Markov chain $(X_t)$,  with the opposite inequality $\tau_0 \leq 1$ being the condition of null recurrence of the chain.  They state (see their Remark 2.6) that ``\textit{it is not known what is the limiting distribution of $\hat{\hat{\mu}}$ when $m_0=1$ and $\tau_0 \leq 1$}". The main difficulty seems to lie on the quantification of the rates of divergence of quantities such as $\sum_{t=1}^n \frac{1}{1+X_{t-1}}$ and $\sum_{t=1}^n \frac{X_t-X_{t-1}-\mu_0}{1+X_{t-1}}$, if they exist. For instance,  \citep{wei1989some} state that if $\tau_0=1$, the convergence rate of $\sum_{t=1}^n \frac{1}{1+X_{t-1}}$ ``\textit{seems to be more difficult; only some upper and lower bounds are given}".\footnote{In the case where $\tau_0<1$, \cite[Theorem 2.18]{wei1989some} show that ${n^{-1+\tau_0}}\sum_{t=1}^n \frac{1}{1+X_{t-1}}$ converges weakly to $c W_{1-\tau_0}$, where $c$ is a constant and random variable $W_{q}$ has the Mittag-Leffler distribution with parameter $q$.  However,   \cite{wei1989some} do not give the value of the constant $c$.  Further,  the asymptotic behavior of $\sum_{t=1}^n \frac{X_t-X_{t-1}-\mu_0}{1+X_{t-1}}$ is also unknown. }
 In summary, the asymptotic theory of the WLS \eqref{1stwls} is, up to now, not fully known, and even if it converges, its convergence rate likely depends on the unknown parameter value $\tau_0$, which makes statistical inference inconvenient.


In this note, we propose an alternative WLS estimator of $\mu$ in model \eqref{gw}, with weights $\frac{1}{t}$ instead of $\frac{1}{1+X_{t-1}}$. The motivation is that the conditional variance of $X_t$ given $X_{t-1}$ is $\sigma_0^2 X_{t-1}+b_0$, and under the unit root assumption, it can be shown that $X_t \sim \mu_0 t$ as $t$ increases to infinity, hence the choice of the weight $\frac{1}{t}$.  Compared to $\sum_{t=1}^n \frac{1}{1+X_{t-1}}$, the harmonic series $ \sum_{t=1}^n  \frac{1}{t}$ is much more tractable since $\sum_{t=1}^n  \frac{1}{t}= \ln n+ \gamma_E+ o(\frac{1}{n})$ as $n$ increases to infinity, where $\gamma_E \approx 0.577$ is the Euler's constant.  We show that the resulting new WLS estimator is consistent for model eq.\eqref{gw}, 
without the additional transience condition eq.\eqref{transient}.  We discuss the implications of this result for statistical inference. 

The rest of the paper is organized as follows. Section 2 proposes the new estimator and derives its asymptotic property for the integrated Galton-Watson process. Section 3 discusses the implications of the new results. Section 4 concludes.
Technical proofs are gathered in the appendix. 
\section{The estimator}
 We consider the $1/t$-weighted WLS estimator:
\begin{equation}
\label{wls1}
(\tilde{m}, \tilde{\mu})':=\arg\min_{(m, \mu)'} \sum_{t=1}^n \frac{1}{t} (X_t-m X_{t-1}-\mu)^2.  
\end{equation}



After centering $\tilde{m}$ and $\tilde{\mu}$ around their respective true values $1$ and $\mu_0$, we get:
\begin{equation}
\label{errorterm}
\begin{bmatrix} \tilde{m}-1  \\
\tilde{\mu}-\mu_0
 \end{bmatrix}
 =\begin{bmatrix}\sum_{t=1}^n \frac{1}{t}X_{t-1}^2  & \sum_{t=1}^n \frac{1}{t} X_{t-1}  \\\sum_{t=1}^n \frac{1}{t} X_{t-1}  & \sum_{t=1}^n \frac{1}{t} \end{bmatrix}^{-1}\begin{bmatrix}
 	\sum_{t=1}^n X_{t-1}W_{t} \frac{1}{t} \\
 	\sum_{t=1}^n W_{t} \frac{1}{t} 
 \end{bmatrix} 
 \end{equation}
where 
 $$W_t= X_{t}-  X_{t-1}- \mu_0$$ is a martingale difference sequence.  We also recall the following well-known scaling property of the integrated Galton-Watson process [see \cite[Theorem 2.1]{wei1989some}]: as $n$ increases to infinity,  the rescaled and time-changed process $(\tilde{X}_{n,s})=(\frac{X_{\lfloor n s \rfloor }}{n}, s>0)$ converges weakly,  in the Skorohod space $\mathcal{D}^{+}([0, \infty[)$ of nonnegative c\`adl\`ag functions on $[0, \infty[$,  towards the Cox-Ingersoll-Ross (CIR) diffusion process:
\begin{equation}
\label{firstdiffusion}
\mathrm{d} Y_s= \mu_0 \mathrm{d} s + \sigma_0 \sqrt{Y_s} \mathrm{d}B_s, \forall s>0, \qquad Y_0=0, 
\end{equation}
where $(B_s)$ is a standard Brownian motion. For our estimation problem, this result can be easily strengthened by introducing the continuous time martingale $(M_{n,s})_s$ :
$$
M_{n,s}= \frac{1}{n}\sum_{t=1}^{\lfloor ns \rfloor } W_t=\frac{X_{\lfloor n s \rfloor }}{n}-\frac{\lfloor n s \rfloor}{n}\mu_0, \qquad \forall s>0.
$$
Then by the continuous mapping theorem (CMT), we get the joint weak convergence of the process:
\begin{equation}
	\label{jointconvergence}
	\Big(\frac{X_{\lfloor n s \rfloor }}{n}, M_{n,s} \Big)\Rightarrow \Big(Y_s, M_s=Y_s- \mu_0s =\sigma_0 \int_0^s \sqrt{Y_s} \mathrm{d}B_s\Big).
\end{equation}


The analysis of the WLS estimator is then based on the following two theorems. 
\begin{thm}
As $n$ increases to infinity, 
\begin{itemize} 
    \item $a)$ $\frac{1}{n} \sum_{t=1}^n  \frac{1}{t}{X_{t-1}}   \Rightarrow \int_0^1 \frac{Y_s}{s} \mathrm{d}s$, where $\Rightarrow$ denotes weak convergence
    \item $b)$ $\frac{1}{n^2} \sum_{t=1}^n \frac{1}{t} {X^2_{t-1}} \Rightarrow \int_0^1 \frac{Y_s^2}{s} \mathrm{d}s,$
\item  $c)$ $\frac{1}{n} \sum_{t=1}^n \frac{1}{t} {W_{t}X_{t-1}} \Rightarrow \sigma_0 \int_{0}^1 \frac{{Y^{3/2}_s}}{s} \mathrm{d}B_s$
\item $d)$ $\frac{1}{\sqrt{\ln n}}  \sum_{t=1}^n \frac{1}{t} {W_{t}}\Rightarrow Z, 
$ where $Z$ follows a normal distribution $N(0,\sigma_0^2 \mu_0)$, and is independent of process $(Y_r)$ 
\item $e)$ The weak convergences in $a)-d)$ are also a joint weak convergence
\end{itemize}
\end{thm}
\noindent We note that the weak limits in $a),b),c)$ all exist, since $\frac{Y_s}{s}$ is continuous on $[0,1]$. Indeed, when $s$ is close to zero, we have 
$
Y_s=\mu_0 s + \sigma_0\int_0^s \sqrt{Y_u} \mathrm{d}u,
$
thus $\frac{Y_s}{s} \rightarrow \mu_0$ almost surely as $s$ decreases to 0. 

The proof of Theorem 1 is quite tedious and involves several limit theorems for stochastic processes. In the following, we provide some intuitions and relegate the formal proof to the appendix. For $a)$, we have:
$$
\frac{1}{n} \sum_{t=1}^n \frac{1}{t} {X_{t-1}} =\int_{1/n}^{1+1/n} \frac{X_{\lfloor ns \rfloor}}{{\lfloor ns \rfloor}} \mathrm{d}s .
$$

  Then the idea is to show that this latter converges to:
$
  \int_0^1 \frac{Y_s}{s} \mathrm{d}s.
$
The difficulty is that $( \frac{X_{\lfloor ns \rfloor}}{n}) \Rightarrow Y_s$ does not directly imply $(\frac{X_{\lfloor ns \rfloor}}{{\lfloor n s \rfloor}}) \Rightarrow Y_s/s$ on the entire domain $(0, \infty)$, due to the singularity at zero. In other words, the CMT cannot be applied directly. Instead, we have to control for the behavior of $ \frac{X_{\lfloor ns \rfloor}}{{\lfloor ns \rfloor}} $ around zero.

$b)$ can be proved in the same way as $a), $ by noting that:
$$
\frac{1}{n^2} \sum_{t=1}^n \frac{1}{t} {X^2_{t-1}} =\int_{1/n}^{1+1/n} \frac{X_{\lfloor ns \rfloor}}{{\lfloor ns \rfloor}  }  \frac{X_{\lfloor ns \rfloor}}{n} \mathrm{d}s.
$$


As for $c)$, we have:
\begin{equation}
\label{limit1}
\frac{1}{n} \sum_{t=1}^n \frac{1}{t}   {W_{t}X_t}   =\int_{1/n}^{1+1/n} \frac{X_{\lfloor ns \rfloor}}{\lfloor ns \rfloor}   dM_{n,s}. 
\end{equation}
This latter is a stochastic integral with respect to the martingale $(M_{n,s})$. The idea is that process $\frac{X_{\lfloor ns \rfloor}}{\lfloor ns \rfloor}$ converges weakly to $(Y_s/s)$, and $(M_{n,s})$ converges weakly to $(M_s)$, and these convergences are in fact a joint convergence. 
This simple idea, however, requires two technical treatments. First, similar as in $a)$, we need to deal with the singularity at 0. Second, \eqref{limit1} is a stochastic integral and for unlike for a deterministic integral, the joint weak convergence of the integrator process and the integrand process is not enough to deduce the weak convergence of the stochastic integral. Instead, some uniformity condition for the martingale $(M_{n,s})$, called Uniform Tightness (UT), needs to be checked.


For $d)$, we write:
\begin{align*}
  \sum_{t=1}^n  {W_{t}}\frac{1}{t}&=  \sum_{t=1}^n \int_{\frac{t}{n}}^{\frac{t+1}{n}}  \frac{1}{\lfloor n s \rfloor }\mathrm{d}M_{n,s}    \approx   \sigma_0 \int_{\frac{1}{n}}^{1+\frac{1}{n}}  \frac{\sqrt{Y_s}}{s} \mathrm{d}B_s.
\end{align*}
This approximation is informal,  because both the bounds and the martingale $(M_{n,s})$ depend on $n$. 
Assuming this is valid, then the stochastic integral $ \sigma_0\int_{a}^1 \frac{\sqrt{Y_s}}{s} \mathrm{d}B_s$ is a local martingale indexed by $a$,  with quadratic variation $\sigma_0^2 \int_{a}^1 \frac{Y_s}{s^2}\mathrm{d}s$.
It is also easily checked that $\frac{\sigma_0^2}{-\ln a}\int_{a}^1 \frac{Y_s}{s^2}\mathrm{d}s$ converges, as $a$ decreases to 0, to $\sigma_0^2\mu_0$, since $Y_s/s \rightarrow \mu_0$ almost surely as $r$ decreases to zero. Hence, by a Central Limit Theorem for stochastic integrals \citep[Theorem 1.19]{kutoyants2013statistical}, we get:
\begin{equation}
\label{limit2}
\frac{\sigma_0}{\sqrt{\ln n}} \int_{1/n}^1 \frac{\sqrt{Y_s}}{s} \mathrm{d}B_s \Rightarrow N(0,\sigma_0^2 \mu_0)
\end{equation}
as $n$ increases to infinity. As a consequence, we also have:
$$
\frac{1}{\sqrt{\ln n}}  \sum_{t=1}^n  {W_{t}} \frac{1}{t} \Rightarrow N(0,\sigma_0^2 \mu_0).
$$
Moreover, as $n$ increases to infinity, $\frac{1}{\sqrt{\ln n}} \int_{1/n}^1 \frac{\sqrt{Y_s}}{s} \mathrm{d}B_s$ depends more and more on the behavior of $(B_s)$ around zero. As a consequence, the joint weak limiting variables in equations \eqref{limit1} and \eqref{limit2} are independent.

Finally,  the fact that the weak convergences $a)$ to $d)$ imply joint weak convergence in $e)$ is a direct consequence of the CMT.  See e.g.  Remark 2.4 of \cite{wei1989some}.

The following theorem is a simple consequence of Theorem 1 and CMT. 
 \begin{thm}
 If there exists a positive exponent $\delta$ such that $\mathbb{E}[Z^{2+\delta}_{i,t}]$ and $ \mathbb{E}[\epsilon_t^{2+\delta}]$ are finite, then the $1/t-$WLS estimator $(\tilde{m},\tilde{\mu}) $ is such that: 
\begin{equation}
\label{mainresult}
 \begin{bmatrix} 
 n(\tilde{m}-1)  \\
\sqrt{\ln n}(\tilde{\mu}-\mu_0)
 \end{bmatrix}
 \Rightarrow
\begin{bmatrix}
 \frac{ \int_0^1  
    \sigma_0 \frac{{Y^{3/2}_s}}{s} \mathrm{d}B_s/}{ \int_0^1 \frac{Y^2_s}{s} \mathrm{d}s }\\
      Z
  \end{bmatrix},
\end{equation}
as $n$ increases to infinity,  where $Z$ follows a normal distribution $N(0, \sigma^2_0 \mu_0)$ and is independent of the diffusion $(Y_s)$. In particular, we have: $\sqrt{\ln n}(\tilde{\mu}-\mu_0) \Rightarrow N(0, \sigma^2_0 \mu_0)$.   
 \end{thm}
 \noindent The convergence rate $(\sqrt{\ln n})$ of our $1/t$-WLS estimator  is the same as \cite{wei1990estimation}'s $\frac{1}{1+X_{t-1}}-$WLS estimator, but the latter is valid only when $\tau_0>1$. In this latter case, \cite[Theorem 2.5]{wei1990estimation} and \cite[Theorem 1.2]{wei1991convergence} show that their WLS has an asymptotic variance of $(\mu_0-\frac{\sigma_0^2}{2})\sigma_0^2$, which is smaller than the asymptotic variance $(\mu_0 \sigma^2_0)$ of our WLS. This means that in case one determines that $\tau_0>1$, then  \cite{wei1990estimation}'s WLS estimator is still more efficient. 
 However, we can always use $\tilde{\mu}$ as the preliminary estimator of $\mu$. Then we can use the consistent estimator $\hat{\sigma}^2$ of $\sigma^2$ in the unit root case suggested by \cite[Theorem 3.5]{winnicki1991estimation}, and estimate $\tau$ by $\hat{\tau}=2\tilde{\mu}/\hat{\sigma}^2$, and check whether it is indeed larger than 1. If this proves the case, then we can use the original estimator of \cite{wei1990estimation} as an improved estimator. If on the other hand $\hat{\tau}$ is non larger than 1, then only our estimator $\tilde{\mu}$ is valid. 
 

  \section{Discussions}
  \subsection{Comparison with the OLS estimator}
As a comparison,   for an integrated Galton-Watson process, the OLS estimator has the following property:
\begin{thm}
The OLS estimator $\hat{\mu}$ is inconsistent,  in the sense that:
 $$
 \begin{bmatrix}n(\hat{m}-1)
 \\
 \hat{\mu}-\mu_0
 \end{bmatrix}
 \Rightarrow \begin{bmatrix}
 \int_0^1 {Y_s^2}  \mathrm{d}s & \int_0^1 {Y_s}  \mathrm{d}s \\
 \int_0^1 {Y_s}  \mathrm{d}s & 1
 \end{bmatrix}^{-1} 
 \begin{bmatrix}
\sigma_0  \int_{0}^1  {Y^{3/2}_s}  \mathrm{d}B_s\\
 \sigma_0 \int_0^1 \sqrt{Y_s} \mathrm{d}B_s
 \end{bmatrix}.
$$
\end{thm}


\subsection{A simulation experiment}
As an illustration,  we consider the model in which $X_t$ is conditionally Poisson given $X_{t-1}$,  with parameter $X_{t-1}+\mu_0$ (hence $\sigma^2=1$).  We set $\mu_0=2$,  and compute the OLS and the $1/t$-WLS estimators on $B=5000$ simulated datasets with sample size $n=100$.  Figure 1 plots the sample histogram of the two estimators.  

\begin{figure}[H]
\centering
\includegraphics[scale=0.15]{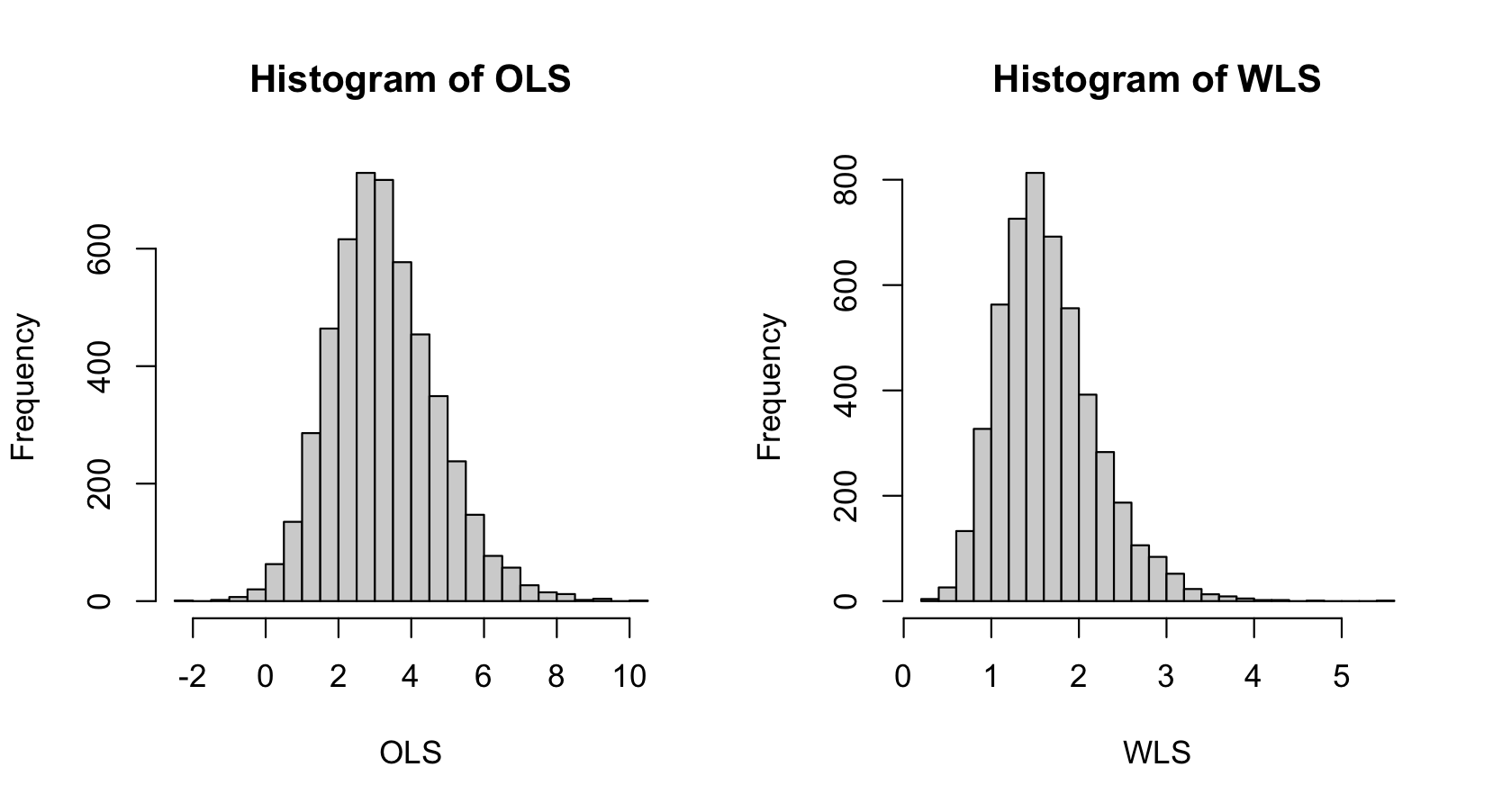}
\caption{Comparison of the sample distributions of the OLS estimator $\hat{\mu}$ (left panel) and the $1/t$-WLS estimator (right panel)}
\end{figure}
We see that the range of the sample distribution of the WLS estimator is significantly narrower compared to its OLS counterpart.  Let us also report some summary statistics.  The OLS estimator has a sample mean of 3.23 (compared to the true value of 2), as well as a sample variance of 2.12,  whereas the WLS estimator has a sample mean of 1.65,  as well as a sample variance of 0.31.  In other words,  despite the slow convergence rate of $\sqrt{\ln n}$,  the WLS estimator is significantly more accurate than the OLS estimator.  If we increase the sample size $n$ to $n=500$ and $n=1000$, respectively,   then the sample mean of the OLS estimator becomes 3.40 and 3.42,  respectively whereas the sample mean of the WLS estimator becomes 1.79 and 1.85,  respectively.  Thus the finite sample bias of the WLS estimator decreases as $n$ increases,  whereas that of the OLS estimator does not.  This confirms the fact that the $1/t$-WLS estimator is consistent, but not the OLS estimator.  

We also report in Table 1 the sample mean of the WLS estimator for different true values of $\mu_0$ and different sample sizes.   
\begin{table}[H]
\centering
\begin{tabular}{|c|c |c |c |c |c |c |c |c  }
\hline
$\mu_0$ & 0.4  & 0.5 & 0.8 & 2 & 3 & 10 \\ \hline
$n=100$ & 0.45 & 0.53 & 0.77 &  1.65 & 2.37  & 7.34  \\ \hline
$n=500$ & 0.45 & 0.53 & 0.79 & 1.79 & 2.55 & 8.10  \\ \hline
\end{tabular}
\caption{Sample mean of the WLS estimator for different $\mu_0$ and $n$. }
\end{table}
The finite sample bias,  normalized by $\mu_0$,  seems to be smaller for small $\mu_0$.  Also, as expected, as $n$ increases,  the relative bias typically diminishes, with the improvement being more significant for large $\mu_0$.  

Finally,  in case the Galton-Watson model is fully parameterized, such as in the INARCH case,  one can also conduct finite sample bias correction,  by first estimating $\mu$,  and then compute the sample bias by simulation as in Table 1.  

\subsection{Uniform inference?}
The lack of unified inference for the Galton-Watson process has been well documented \citep{heyde1974notes, wei1990estimation}.  Can we also use the new WLS estimator proposed in this paper in the stationary case ($m<1$)? We have the following theorem:
\begin{thm}
If $m<1$,  then we have:
\begin{itemize} 
    \item $a)$ $\frac{1}{\ln n } \sum_{t=1}^n  \frac{1}{t}{X_{t-1}}   \rightarrow \lim_{t \to \infty } \mathbb{E}[X_t]= \frac{\mu_0}{1-m_0}$,  almost surely, 
    \item $b)$ $\frac{1}{\ln n } \sum_{t=1}^n \frac{1}{t} {X^2_{t-1}} \rightarrow \lim_{t \to \infty } \mathbb{E}[X_t^2],$  almost surely, 
\item  $c)$ $\sum_{t=1}^n \frac{1}{t} {W_{t}X_t} $ and $ \sum_{t=1}^n \frac{1}{t} {W_{t}}$ both converge in $L^2$ to some random variables $Z_1$,  $Z_2$.  In particular,  they also converge jointly and weakly to $(Z_1,Z_2)$. 
\end{itemize}
and $$   \begin{bmatrix}
(\ln n)(\tilde{m}-m_0 )\\
(\ln n)(\tilde{\mu}-\mu_0)
\end{bmatrix} \Rightarrow 
\begin{bmatrix}
\lim_{t \to \infty } \mathbb{E}[X^2_t] & \lim_{t \to \infty } \mathbb{E}[X_t] \\
\lim_{t \to \infty } \mathbb{E}[X_t] & 1
\end{bmatrix}^{-1}
\begin{bmatrix}
Z_1 \\
Z_2
\end{bmatrix},
$$
where $\lim_{t \to \infty } \mathbb{E}[X_t^2]$ is the second order moment of the stationary distribution of $(X_t)$ and is given by:
$
\mathbb{E}[X_t^2]= \Big(\frac{\mu_0}{1-m_0} \Big)^2+\frac{\sigma_0^2\frac{\mu_0}{1-m_0}+b_0}{1-m_0^2}.
$
\end{thm}
This theorem is straightforward and the proof is omitted.  

Thus the new WLS estimator is consistent in both the stationary and  integrated case,  but its convergence rate is sub-optimal in the stationary case compared to other estimators such as the OLS estimator\footnote{Or the maximum likelihood (ML) estimator, if the conditional distribution is fully specified.} in the case where $m_0<1$,  which converges at the standard parametric rate $\sqrt{n}$.  As a consequence,  the OLS is still preferred,  if the researcher concludes that the process $(X_t)$ is stationary.  
\subsection{Stationarity and unit root tests}
Section 3.3 leads to a very natural question: In practice, how to tell from the data whether process $(X_t)$ is stationary or integrated? One can use the following procedure,  first advocated by \cite{kwiatkowski1992testing} (KPSS): 

$i)$ Conduct a stationarity test,  such as the KPSS test,  with null hypothesis $m_0<1$ and alternative hypothesis $m_0=1$.  This test has for instance been suggested by \cite{michel} for count-valued processes. 

$ii)$ Complement the stationarity test with a unit root test on $m$, using the limit theory of Theorem 2 or Theorem 3 for $\hat{m}$ or $\tilde{m}$.  In both cases, the limiting distribution involves parameters $\mu$ and $\sigma^2$,  which have to be estimated.  The new $1/t-$WLS allows to estimate $\mu$, whereas $\sigma^2$ can be estimated using \cite[Theorem 3.5]{winnicki1991estimation}. If the unit root test is not rejected,  but the stationarity test is rejected, then the process is likely nonstationary.  

$iii)$ If the unit root test is rejected,  but not the stationarity test,  then the process is likely stationary.  In this case,  the OLS estimator should be used to estimate $m$ and $\mu$,  since it converges much faster than the WLS estimator in the stationary case.  See also  \cite[Proposition 3.1]{winnicki1991estimation} for an estimator of $\sigma^2$ and $b$ and its properties in this stationary case.  


\section{Conclusion}
This paper has introduced a new estimator of the intercept parameter $\mu$ of an integrated Galton-Watson process, that is consistent regardless of whether $\tau_0=2\mu_0/\sigma_0^2$ is larger than 1.  The weight function $1/t$ proposed in this paper is new in this literature, and it can be potentially very useful, 
since similar inconsistency issues exist also in many more complicated models. 
For instance,  \cite{barreto-souza2023} consider an INARCH model in which parameter $m_n$ is local-to-unity (\textit{i.e.} $m_n=1+\frac{\gamma}{n}$),  which extends the exactly integrated INARCH model by allowing $\gamma$ to be potentially non-zero.  Their analysis assumes $\mu$ to be ex ante fixed, likely to avoid the inconsistency issue. 
 \cite{barczy2014asymptotic} consider the INAR(2) model:\begin{equation}
 	\label{inarp}
X_{-1}=X_0=0, \qquad  X_t= \alpha_1 \circ X_{t-1} +   \alpha_2  \circ X_{t-2} + \epsilon_t, \forall t,\end{equation} where probability parameters $\alpha_1, \alpha_2$ are both positive, the count sequence $(\epsilon_t)$ is i.i.d. across $t$,  with a mean $\mu$. \cite[Theorem 2.1]{barczy2014asymptotic} show that under the unit root assumption $ \alpha_1+   \alpha_2 =1,$ the OLS estimator of $\mu$ is inconsistent. Their model can be viewed as a second-order extension of \eqref{gw}. \cite{michel2020limiting} derived the scaling property towards the CIR diffusion of the non-Markov INGARCH(1,1) process, but stop short of deriving a consistent estimator for $\mu$,  among others. This model can be viewed as the infinite order extension of \eqref{gw}. The analyses of these models are left for future research.   

\bibliographystyle{apalike}
\bibliography{lib}
\appendix

\newpage

\setcounter{page}{1}
\section{Proof of Theorem 1}
For $a)$. We remark that:
\begin{itemize}
	\item Since $\mathbb{E}[X_t]=t \mu_0$, we get, for any $n$, and any given $\eta>0$:
	$$\mathbb{E}[\frac{1}{n} \sum_{t=1}^{\lfloor n \eta \rfloor } \frac{1}{t} {X_{t-1}} ]    \leq \frac{{\lfloor n \eta \rfloor +1}}{n}   \mu_0 \leq 2 \eta \mu_0$$ 
	Thus $\frac{1}{n} \sum_{t=1}^{\lfloor n \eta \rfloor } \frac{1}{t} {X_{t-1}} $ converges to zero in probability uniformly in $n$, as $\eta$ decreases to zero. 
	\item By the CMT, for a fixed $\eta>0$, we have $\frac{1}{n} \sum_{t=\lfloor n \eta \rfloor }^{n} \frac{1}{t} {X_{t-1}}\Rightarrow   \int_{\eta}^1 \frac{Y_s}{s} \mathrm{d}s$ 
	\item finally, $\int_0^{\eta} \frac{Y_s}{s} \mathrm{d}s \rightarrow 0$ almost surely as $\eta$ goes to 0. 
\end{itemize}
Thus we conclude that $$
\frac{1}{n} \sum_{t=1}^n \frac{1}{t} {X_{t-1}}  \Rightarrow   \int_0^1 \frac{Y_s}{s} \mathrm{d}s
$$
as $n$ goes to infinity. Hence $a)$ holds. 

For $c)$, let us simply show that, for any given $\epsilon>0$, we have  
\begin{equation}
	\label{weakconvergenceUT}
\int_{\epsilon}^{1} \frac{X_{\lfloor ns \rfloor}}{\lfloor ns \rfloor}   dM_{n,s} \Rightarrow \int_\epsilon^1 \frac{Y_s}{s} \mathrm{d}M_s
\end{equation}
 By \citep[Theorem 2.6]{jakubowski1989convergence} [see also \cite{kurtz1991weak} and \cite[Chapter VI, section 6]{jacod2013limit}], we need to check that the sequence of martingales $(M_{n,s}, \epsilon<s<1)$ satisfy the UT condition. Since $(W_t)$ is a $(M_{n,s})$ is a sequence of martingales,  by \cite[Proposition 3.2]{jakubowski1989convergence}, it suffices to check that:
$$
\sup_n \mathbb{E}[\displaystyle \max_{1 \leq t \leq n}  \frac{|W_t|}{n}] < \infty. 
$$
We remark $\max_{1 \leq t \leq n}  \frac{|W_t|}{n} \leq \max_{1 \leq t \leq n} \frac{|X_t|}{n}+ \frac{\mu_0}{n}$. So it suffices to prove that $\sup_n \mathbb{E}[\displaystyle \max_{1 \leq t \leq n}  \frac{|X_t|}{n}] < \infty. $ Since $(X_t)$ is a sub-martingale, by Doob's $L^q$ inequality,  we have:
$$
\mathbb{E}[  (\max_{1 \leq t \leq n}  {|X_t|} )^q] \leq (\frac{q}{q-1})^q \mathbb{E}[X_n^q] 
$$
for any exponent $q>1$.  But by the weak convergence of the process $(\frac{X_{\lfloor ns \rfloor }}{n})$ towards $(Y_s)$, we have: $\mathbb{E}[X_n^q]=O(n^q)$ for any $q$ between 1 and $1+\delta$.  Thus 
the sequence $ \mathbb{E}[\max_{1 \leq t \leq n}  \frac{|W_t|}{n}] $,  $n$ varying,  is bounded. Thus the UT condition is satisfied and \eqref{weakconvergenceUT} holds.

For $d)$, let us apply the martingale central limit theorem to the triangular array $W_{t,n}= \frac{1}{\sqrt{\ln n }} W_{t}  \frac{1}{t} $. We check that:
\begin{itemize} 
    \item $i)$  $V_n:=\sum_{t=1}^n \mathbb{V}ar[W_{t,n} \mid X_{t-1},X_{t-2}, \cdots  ] \rightarrow \mu_0$ in probability as $n$ increases to infinity. 
    \item $ii)$  (Lindeberg condition): $\sum_{t=1}^n \mathbb{E}[W^2_{t,n} \mathbbm{1}_{X_{t,n}>\epsilon \sqrt{V_n} } \mid X_{t-1},X_{t-2}, \cdots] \rightarrow 0$ in probability for any given $\epsilon>0$. 
\end{itemize}
Indeed, if these two conditions are satisfied, then  $\frac{1}{\mu_0}\sum_{t=1}^n W_{t,n}$ converges weakly to a  normal variable. 

For $i)$, we have:
\begin{align*}
    \sum_{t=1}^n \mathbb{V}ar[X_{t,n} \mid X_{t-1} ,X_{t-2}, \cdots]&= \frac{1}{\ln n}   \sum_{t=1}^n    \frac{\sigma_0^2 X_{t-1}+b_0}{t^2}\\
&= \frac{1}{\ln n}   \sum_{t=1}^n    \frac{\sigma_0^2 \mathbb{E}[X_{t-1}]+\sigma_0^2(X_{t-1}-\mathbb{E}[X_{t-1}])+b_0}{t^2}\\
    &= \frac{1}{\ln n}   b_0 \sum_{t=1}^n  \frac{1}{t^2} +  \frac{1}{\ln n}   \mu_0 \sum_{t=1}^n  \frac{1}{t}  +   \frac{1}{\ln n}\sigma_0^2   \sum_{t=1}^n  \frac{X_{t-1}-(t-1)\mu_0 }{t^2},
\end{align*}
where we have used $\mathbb{E}[X_{t-1}]=\mu_0( t-1)$.  The first and second terms converge to 0 and $\mu_0$, respectively, as $n$ increases to infinity. For the second term, we can show that it converges to zero in $L^2$, and thus also in probability. Thus $\sum_{t=1}^n \mathbb{V}ar[W_{t,n}\mid X_{t-1},X_{t-2}, \cdots] \rightarrow b_0$ in probability. 

For $ii)$,  by Markov's inequality,  we have, 
\begin{align*}
    \sum_{t=1}^n \mathbb{E}[|W_{t,n}|^{2+\delta} \mathbbm{1}_{W_{t,n}>\epsilon \sqrt{V_n} } \mid X_{t-1},X_{t-2}, \cdots] &=   \sum_{t=1}^n \mathbb{E}[|W_{t,n}|^{2+\delta} \mathbbm{1}_{|W_{t,n}||^{\delta}>\epsilon^{\delta} {V_n}^{\delta/2} } \mid X_{t-1},X_{t-2}, \cdots] \\
     &\leq   \frac{1}{\ln n} \sum_{t=1}^n  \frac{1}{t^2} \frac{\mathbb{E}[|W_t|^{2+\delta}]}{\epsilon^{\delta} (\ln n)^{\delta/2} V_n^{\delta/2} t^{\delta}}.
    \end{align*}
Thus it suffices to show that $\sum_{t=1}^n \mathbb{E}[|W_t|^{2+\delta}]/t^{2+\delta}$ is bounded in order for the right hand side (RHS) of the above inequality to go to zero as $n$ increases to infinity. To this end, let us show that:
\begin{equation}
	\label{target}
 \mathbb{E}[|W_t|^{2+\delta}]=O(t^{1+\delta/2}).
 \end{equation} We write: $$W_t=\underbrace{\sum_{i=1}^{X_{t-1}} (Z_{i,t}-1)}_{:=A_t}+ \underbrace{\epsilon_t-\mu_0}_{:=B_t}.$$ By Jensen's inequality, we have: $
|W_t|^{2+\delta} \leq 2^{1+\delta} ( A_t^{2+\delta}+ B_t^{2+\delta}).
$
By taking expectation, we get:
$$
\mathbb{E}[|W_t|^{2+\delta}] \leq 2^{1+\delta} \mathbb{E}[A_t^{2+\delta}]+ 2^{1+\delta} \mathbb{E}[B_t^{2+\delta}].
$$
The second term on the RHS is a constant. For the first term, we use the Rosenthal inequality, we deduce that there exist a positive constant $C_1$ such that:
$$
\mathbb{E}[\left|A_t\right|^{2+\delta}|X_{t-1}]
\;\le\;
C_1\left[
X_{t-1}\,\mathbb{E}\left|Z-1\right|^{2+\delta}
\;+\;
\left(X_{t-1}\,\mathrm{Var}(Z)\right)^{(2+\delta)/2}
\right].
$$
Thus by taking marginal expectations, and by using the fact that $x^{(2+\delta)/2}+1 \geq x$ for any $x>0$, we get:
$$
\mathbb{E}[\left|A_t\right|^{2+\delta}] \leq C_2 \mathbb{E}[X^{(2+\delta)/2}_{t-1}]+C_3,
$$
where $C_2$, $C_3$ are positive constants.  Then by the weak convergence of $\frac{X_{\lfloor ns \rfloor }}{s}$, we can show that:
$
\mathbb{E}[X_{t-1}^{1+\delta/2}]=O(t^{1+\delta/2}).
$ Thus \eqref{target} is satisfied. The result follows.
    \section{Proof of Theorem 2}
     Let us define $F_n$ and $d_n$ by:
     $$
    F_n= \begin{bmatrix}\sum_{t=1}^n \frac{1}{t}X_{t-1}^2  & \sum_{t=1}^n \frac{1}{t} X_{t-1}  \\\sum_{t=1}^n \frac{1}{t} X_{t-1}  & \sum_{t=1}^n \frac{1}{t} \end{bmatrix}, \qquad 
    d_n=\begin{bmatrix}
     	\sum_{t=1}^n X_{t-1}W_{t} \frac{1}{t} \\
     	\sum_{t=1}^n W_{t} \frac{1}{t} 
     \end{bmatrix},
 $$
     Then we rescale them and define $\tilde{F}_n$ and $\tilde{d}_n$ through:
  $$
  \tilde{F}_n= \begin{bmatrix}
  n^{-1} & 0 \\
  0 & 1/\sqrt{ \ln n}
  \end{bmatrix} F_n  \begin{bmatrix}
  n^{-1} & 0 \\
  0 & 1/\sqrt{  \ln n}
  \end{bmatrix}= \begin{bmatrix}
  n^{-2} \sum_{t=1}^n X_{k-1}^2  \frac{1}{t}  & \frac{1}{n \sqrt{\ln n}} \sum_{t=1}^n X_{t-1} \frac{1}{t} \\
  \frac{1}{n \sqrt{\ln n}} \sum_{t=1}^n X_{t-1} \frac{1}{t} & \frac{1}{\ln n}\sum_{t=1}^n  \frac{1}{t} 
  \end{bmatrix},
  $$
  and
  $$
  \tilde{d}_n=\begin{bmatrix}
  n^{-1} & 0 \\
  0 & 1/\sqrt{ \ln n}
  \end{bmatrix} d_n=\begin{bmatrix}
  n^{-1}  \sum_{t=1}^n X_{t-1}W_{t} \frac{1}{t}    \\
   \frac{1}{\sqrt{ \ln n}}  \sum_{t=1}^n W_{t} \frac{1}{t} 
  \end{bmatrix}.
  $$ 
  Then by Theorem 1, the joint distribution of $(\tilde{F}_n, \tilde{d}_n)$ converges weakly to the distribution of $(\tilde{F},\tilde{d})$ given by:
  $$
 \tilde{F}= \begin{bmatrix}
  \int_0^1 \frac{Y^2_s}{s} \mathrm{d}s&0 \\
0 & 1
  \end{bmatrix}, \qquad \tilde{d}=\begin{bmatrix}
  \int_0^1  
    \sigma_0 \frac{{Y^{3/2}_s}}{s} \mathrm{d}B_s \\
      Z
  \end{bmatrix}.
  $$
  Thus we get:
  $$
  F_n^{-1} d_n= \begin{bmatrix}
  n^{-1} & 0 \\
  0 & 1/\sqrt{ \ln n}
  \end{bmatrix} \tilde{F}^{-1}_n \begin{bmatrix}
  n^{-1} & 0 \\
  0 & 1/\sqrt{ \ln n}
  \end{bmatrix} \begin{bmatrix}
  n & 0 \\
  0 & \sqrt{  \ln n}
  \end{bmatrix} \tilde{d}_n,
  $$
  or $$
  \begin{bmatrix}
  n & 0 \\
  0 & \sqrt{ \ln n}
  \end{bmatrix}   F_n^{-1} d_n = \tilde{F}^{-1}_n\tilde{d}_n  \Rightarrow \tilde{F}^{-1}\tilde{d}
  $$
  by the CMT.   In particular, $\sqrt{\ln n}(\tilde{\mu}-\mu)$ converges to a normal distribution. 
%

 \section{Proof of Theorem 3}
The OLS estimator satisfies:
\begin{equation}
	\label{error2}
	\begin{bmatrix}
		n(\hat{m}-1)
		\\
		\hat{\mu}-\mu_0
	\end{bmatrix}
	= \begin{bmatrix}\sum_{t=1}^n X_{t-1}^2  & \sum_{t=1}^n   X_{t-1}  \\\sum_{t=1}^n  X_{t-1}  & n  \end{bmatrix}^{-1}
	\begin{bmatrix}
		\sum_{t=1}^n X_{t-1}W_{t}   \\
		\sum_{t=1}^n W_{t} 
	\end{bmatrix}.
\end{equation}

Then it suffices to remark the following weak convergences:
\begin{itemize} 
	\item $\frac{1}{n^2} \sum_{t=1}^n  {X_{t-1}}   \Rightarrow \int_0^1 {Y_s} \mathrm{d}s$
	\item $\frac{1}{n^3} \sum_{t=1}^n  {X^2_{t-1}} \Rightarrow \int_0^1 {Y_s^2}  \mathrm{d}s,$
	\item  $\frac{1}{n^2} \sum_{t=1}^n  {W_{t}X_t} \Rightarrow \sigma_0 \int_{0}^1  {Y^{3/2}_s}  \mathrm{d}B_s$
	\item  $\frac{1}{n}  \sum_{t=1}^n  {W_{t}}\Rightarrow  \sigma_0 \int_0^1 \sqrt{Y_s} \mathrm{d}B_s$.
\end{itemize}
Here, the first two convergences are the analogues of properties $a)$ and $b)$ in Theorem 1 and are immediate consequences of the CMT. The next two convergences are similar to $c)$ of Theorem 1.  Finally, by the CMT,  the above weak convergences is a joint convergence.    The result follows.

\end{document}